\documentclass[12pt, a4paper]{article}
\usepackage{custom_article}

\author{Andriy Haydys\\
}
\title{$\rG_2$ instantons and the Seiberg--Witten monopoles}
\date{18th March 2017}

\begin{document}
\maketitle

%===============================================================================
 \begin{abstract}
  I describe a relation (mostly conjectural) between the Seiberg--Witten monopoles, Fueter sections, and $\rG_2$ instantons. 
  In the last part of this article I gathered some open questions connected with this relation.
 \end{abstract}
%===============================================================================

%===============================================================================
\section{Introduction}
%===============================================================================

\pagestyle{fancy}
\fancyhead[LE]{}
\fancyhead[CE]{$\rG_2$ instantons and the Seiberg--Witten monopoles}
\fancyhead[RE]{}
\fancyhead[LO]{}
\fancyhead[CO]{$\rG_2$ instantons and the Seiberg--Witten monopoles}
\fancyhead[RO]{}

By a celebrated result of Berger~\cite{Berger55_HolGps}, the holonomy group of a simply connected irreducible nonsymmetric Riemannian manifold is one of the following:
\[
\SO(n),\ \U(n),\ \SU(n),\ \Sp(n),\ \Sp(n)\Sp(1),\ \rG_2,\ \text{and }\; \Spin(7).
\]  
Manifolds with holonomies in $\SU(n),\ \Sp(n),\ \rG_2,\ \text{or } \Spin(7)$ are distinguished by the fact that these admit non-trivial covariantly constant spinors~\cite{Wang89_ParallelSpinors}. 
This is intimately related to the concept of supersymmetry and therefore manifolds with these holonomy groups are particularly important in physics.    
Manifolds with holonomies in $\SU(n)$ and $\Sp(n)$ are K\"ahler and therefore amenable to methods of complex geometry. 
The last two groups --- the exceptional holonomy groups --- are related to the algebra of octonions and can occur only in dimensions $7$ and $8$ respectively. 
In this article I focus on $\rG_2$ manifolds although much of what appears below has analogues in the setting of $\Spin(7)$ manifolds.
A reader interested in special holonomy and in particular exceptional groups is encouraged to consult~\cite{Salamon03_TourExceptGeom, Bryant87_MetrWithExceptHol, Joyce07_RiemannianHolGpsCalibr, Salamon:89}.

The existence of manifolds with exceptional holonomies was established much later after Berger's discovery
 \cites{Bryant87_MetrWithExceptHol, BryantSalamon89_CompleteExcept, Joyce96_CompactG2I_II, Joyce96_CompactSpin7}. 
The twisted connected sum construction~\cite{Kovalev03_TCS} and its refinement~\cite{CHNP15_G2Mflds} turned out to be very efficient in  constructing compact $\rG_2$ manifolds. 
Moreover, it is likely that there are many $7$-manifolds with disconnected moduli space of $\rG_2$ metrics.
At present, only a few examples are known~\cite{CGN_AnalytInvOfG2mflds_arx}.
Higher dimensional gauge theory pioneered in~\cite{DonaldsonThomas:98} may be a useful tool in studies of the moduli spaces of metrics with exceptional holonomies.

\medskip 

Although one can define $\rG_2$ as the automorphism group of the octonion algebra, in the setting of Riemannian geometry another definition is usually more convenient.
Namely, we define $\rG_2$ as the stabilizer group in $\GL(7;\R)$ of a generic 3-form, which can be chosen for example as
\[
\varphi_0=\sum_{i=1}^3\om_i^+\wedge dy_i + dy_1\wedge dy_2\wedge dy_3.
\] 
Here $\R^7$ is thought of as $\R^4\oplus \R^3$, $(\om_1^+,\om_2^+,\om_3^+)$ is the standard basis of self-dual 2-forms on $\R^4$, and $(y_1, y_2, y_3)$ are coordinates on $\R^3$.
It turns out that in fact $\rG_2\subset \SO(7)$, in particular $\rG_2$ is compact. 
Since the orbit of $\varphi_0$ is open in $\Lambda^3(\R^7)^*$, we have $\dim\rG_2=\dim\GL(7;\R)-\dim\Lambda^3(\R^7)^*=14$.
These facts can be found for instance in~\cite{Bryant87_MetrWithExceptHol,BryantSalamon89_CompleteExcept}.

Given any subgroup $G\subset\SO(n)$ we obtain a splitting of $\mathrm{Lie}(SO(n))=\mathfrak{so}(n)$ as a direct sum of vector spaces: $\mathfrak g\oplus \mathfrak g^\perp$. 
Here $\mathfrak g$ is the Lie algebra of $G$. 
This means that if the holonomy group of a Riemannian $n$-manifold reduces to a subgroup, then the bundle of two forms splits.
In the case of a seven-manifold $Y$ with holonomy group $\rG_2$ we have $\Lambda^2T^*Y=\Lambda^2_{14}T^*Y\oplus \Lambda^2_7\,T^*Y$, where the subscript indicates the rank of the corresponding subbundle. 
Denote by $\pi_7\colon \Lambda^2T^*Y\to \Lambda^2_7\,T^*Y$ the corresponding projection.

Let $A$ be a connection on a principal bundle over $Y$ with the curvature $F_A$. 
Then $A$ is called a $\rG_2$ instanton, if  $\pi_7(F_A)=0$.
It turns out   that  $\rG_2$ instantons have a lot in common with flat connections over three-manifolds. 
For example, both $\rG_2$ instantons and flat connections over three-manifolds are critical points of Chern--Simons-type functionals.
Following this analogy it is tempting to define a Casson-type invariant of compact $\rG_2$ manifolds by counting $\rG_2$ instantons in a suitable sense.   
This idea has been suggested in~\cites{DonaldsonThomas:98} and later some details were worked out in~\cite{DonaldsonSegal:09}.

One of the central problems of the above approach is the non-com\-pact\-ness of the moduli space of $\rG_2$ instantons.
Even though a substantial progress in understanding the Uhlenbeck-type compactification has been made~\cites{Price83_MonotonicityFormula, Nakajima88_CompHighDim, Tian:00}, some fundamental problems remain open.
  
It is expected that a compactification of the moduli space of $\rG_2$ instantons can be constructed by adding certain ``ideal instantons''.
One of the ingredients of ideal instantons~\cites{DonaldsonSegal:09,Haydys:12_GaugeTheory_jlms} are Fueter sections, which are closely related to a generalization of the \SWn equations~\cites{Taubes:99}. 
Hence, one may expect that there could be a relation between $\rG_2$ instantons and the Seiberg--Witten monopoles.

Let me explain briefly  what such a relation might be useful for (somewhat more details are given at the end of Section~\ref{Sect_RelationG2SW}). 
As it was observed in~\cite{DonaldsonSegal:09}, the na\"{i}ve count of $\rG_2$ instantons on a compact $\rG_2$ manifold can not produce a deformation-invariant number 
but rather this number will jump in a finite number of points as one changes the $\rG_2$ metric in a one parameter family.
The jumps are closely related to degenerations of $\rG_2$ instantons to  Fueter sections supported on  certain associative submanifolds. 
Since Fueter sections also appear as degenerations of the Seiberg--Witten monopoles, one may hope that by counting $\rG_2$ instantons together with the Seiberg--Witten monopoles on associative submanifolds the mutual degenerations will cancel in pairs. 
How close this picture is to the reality, if at all, remains for the future studies.

\medskip

This article is organized as follows. 
In Sections~\ref{Sect_FuetMapsandSW} and~\ref{Sect_CompactnessSW}  I describe briefly a relation between the Seiberg--Witten equations and Fueter sections in dimension three. 
In Section~\ref{Sect_CompG2Instant} I give some details on the relation between $\rG_2$ instantons and Fueter sections.
A conjectural qualitative picture of a relation between $\rG_2$ instantons and the Seiberg--Witten monopoles is outlined in Section~\ref{Sect_RelationG2SW}.
The last section of this article is devoted to some open questions related to the material of the preceding sections.  

\medskip

\textsc{Acknowledgements.} I am thankful to A.~Doan for reading the draft of this manuscript and for helpful suggestions.

\section{Fueter maps and generalized Seiberg--Witten equations}\label{Sect_FuetMapsandSW}

To keep the exposition clear, in this section I restrict myself to the case of $\R^3$ as a background three-manifold. 
All constructions and statements appearing below can be generalized to arbitrary oriented Riemannnian three-manifolds. 
Moreover, there are also corresponding generalizations for oriented Riemannian four-manifolds.
Some details can be found for example in~\cite{Haydys15_DiracOperatorsInGauge}.

\medskip   

Let $(\mathcal M, I_1, I_2, I_3)$ be an almost hypercomplex manifold, i.e., $(I_1, I_2, I_3)$ is a triple of complex structures satisfying the quaternionic relations.
We say that a map $u\colon \mathbb R^3\to \mathcal M$ is Fueter if 
\begin{equation}\label{Eq_FueterR3}
I_1\frac{\partial u}{\partial x_1} + I_2\frac{\partial u}{\partial x_2} + I_3\frac{\partial u}{\partial x_3}=0,
\end{equation}
where $(x_1, x_2, x_3)$ are coordinates on $\mathbb R^3$. 
For example, if $\mathcal M =\mathbb H$, then $u$ is a harmonic spinor on $\R^3$. 

To obtain a gauged version of~\eqref{Eq_FueterR3}, we require that $\mathcal M$ is hyperK\"ahler. 
Assume also that a Lie group $\rG$ acts on $\mathcal M$ preserving its hyperK\"ahler structure and denote by $K_\xi$ the Killing vector field corresponding to $\xi\in\fg=\mathrm{Lie}(\rG)$: $K_\xi(m)=\left.\tfrac d{dt}\right |_{t=0} \exp(t\xi)\cdot m$.  
Moreover, we assume that there is a hyperK\"ahler moment map $\mu\colon\cM\to\fg\otimes\mathbb R^3$, where $\fg$ is identified with $\fg^*$ via a fixed $\mathrm{ad}$-invariant scalar product.

Given $a\in\Omega^1(\mathbb R^3;\fg)$, which should be thought of as a connection on the product $\rG$-bundle, and a map $u\in C^\infty(\mathbb R^3; \cM)$ we define its covariant derivative by $\nabla^a u= du + K_{a(\cdot)}(u)$.  
Then the generalized Seiberg--Witten equations corresponding to $(\cM, \rG)$ read:
\begin{equation}\label{Eq_gSW}
 I_1\nabla_1^a u + I_2\nabla_2^au +I_3\nabla_3^a u=0 \quad\text{and}\quad F_a=\mu\comp u.
\end{equation}
Notice that in the second equation the identification $\Omega^2(\mathbb R^3;\fg)\cong C^\infty(\mathbb R^3; \fg\otimes\mathbb R^3)$ is used.

For example, the choice $(\cM, \rG) =(\H,\U(1))$ yields the classical Seiberg--Witten equations. 
Choosing $(\cM, \rG)=(\mathbb H^n, \U(1))$, where $\U(1)$ acts componentwise, one clearly obtains a generalization of the classical Seiberg--Witten equations involving $n$ spinors instead of $1$.
These equations are studied in~\cite{HaydysWalpuski15_CompThm_GAFA} and are also discussed in the next section in some detail. 
Generalized \SWn equations with an infinite dimensional target space are studied in~\cite{Haydys:12_GaugeTheory_jlms}, see also Section~\ref{Sect_CompG2Instant}.
Other examples can be found in~\cite{Haydys15_DiracOperatorsInGauge}.

\medskip

A link between solutions of generalized Seiberg--Witten equations and Fueter maps can be constructed as follows~\cite{Haydys:12_GaugeTheory_jlms}.
Let $(a_\e, u_\e)$ be a solution of
\begin{equation} \label{Eq_gSW_e}
 I_1\nabla_1^a u_\e + I_2\nabla_2^au_\e +I_3\nabla_3^a u_\e=0 \quad\text{and}\quad \e^2F_{a_\e}=\mu\comp u_\e,
\end{equation}
where $\e\in (0,1]$ is a parameter. 
Assume that there exists a sequence $\e_i$ converging to $0$ such that the corresponding sequence of solutions $(a_{\e_i}, u_{\e_i})$ converges to some limit $(a_0, u_0)$ in an appropriate functional space. 
Then we obtain 
\begin{equation} \label{Eq_gSW_0}
 I_1\nabla_1^a u_0 + I_2\nabla_2^au_0 +I_3\nabla_3^a u_0=0 \quad\text{and}\quad \mu\comp u_0=0. 
\end{equation}

Suppose that $0\in\fg\otimes\mathbb R^3$ is a regular value of the momentum map $\mu$ and
$\rG$ acts freely on $\mu^{-1}(0)$.
Then $\cM_0=\mu^{-1}(0)/\rG$ inherits a hyperK\"ahler structure.
Furthermore, for any solution $(u_0, a)$ of~\eqref{Eq_gSW_0} denote by $v_0\colon\R^3\to \cM_0$ the projection of $u_0$.
Then the following holds.
\begin{thm}[\cite{Haydys:12_GaugeTheory_jlms}*{Thm~4.2}]\label{Thm_FuetMapsHKred}
 The map $(u_0, a)\mapsto v_0$ yields a  bijective correspondence between the moduli space
of solutions of~\eqref{Eq_gSW_0} and the space of Fueter maps $\mathbb R^3\to \cM_0$.\qed
\end{thm}

This result provides a basis for links between different gauge-theoretic problems. Some details are given in Section~\ref{Sect_CompG2Instant}.
In some sense, this is an analogue of the relation between symplectic vortices and pseudoholomorphic curves~\cite{SalamonEtAl:02}*{Sect.\,2.4}.

To explain a mechanism leading to the appearance of the small parameter $\e$, assume that the target space $\cM$ is a quaternion-Hermitian \emph{vector space} and the moment map is homogeneous of degree two. 
Suppose there is a sequence of solutions $(a_i, u_i)$ such that the sequence $r_i:=\| u_i\|_{L^2}$ diverges to $+\infty$.
If we set $(a_i, \tilde u_i)=(a_i, r_i^{-1}u_i)$, then 
\[
 I_1\nabla_1^{a_i} \tilde u_i + I_2\nabla_2^{a_i} \tilde u_i +I_3\nabla_3^{a_i} \tilde u_i=0 \quad\text{and}\quad 
 r_i^{-2}F_{a_i}=\mu\comp \tilde u_i,
\] 
i.e., $(a_i, \tilde u_i)$ solves~\eqref{Eq_gSW_e} with $\e=r_i^{-1}$.
This is intimately related to the compactness problem for the moduli space of generalized \SWn monopoles as it is described in the next section. 

\medskip

Consider the following example: $\rG$ is an arbitrary compact Lie group and  $\cM=\fg\otimes\H$. 
A simple computation shows that in this case solutions of~\eqref{Eq_gSW} can be interpreted as flat stable $\rG_\C$--connections~\cite{Haydys15_DiracOperatorsInGauge}.
Here a $\rG_\C$--connection $\cA=a+ ib$ is called stable, if $d_a^* b=0$.    
In particular, for $\rG=\SO(3)$ we obtain flat stable $\PSL(2,\C)$--connections.
The compactness property of the corresponding moduli space was studied in~\cite{Taubes13_PSL2Ccmpt, Taubes15_CorrigendumToPSL2C}. 
Roughly speaking, for this particular choice of $(\cM, \rG)$ Taubes shows that any sequence $(a_i, u_i, \e_i)$ of solutions of~\eqref{Eq_gSW_e} with $\e_i\to 0$ and $\| u_i\|_{L^2}=1$ has a convergent subsequence in a suitable sense.
Hence, it is not that surprising that the methods pioneered in~\cite{Taubes13_PSL2Ccmpt} are also applicable in other instances of generalized Seiberg--Witten equations.
This has been already confirmed in a number of cases~\cite{Taubes13_CxASD_Arx, HaydysWalpuski15_CompThm_GAFA, Taubes16_SWDim4_Arx,  Taubes17_VafaWitten_Arx}.

\begin{remark}
Equations~\eqref{Eq_nSW0} make sense on any three--manifold provided $\cM$ admits a permuting action of $\Sp(1)$~\cite{Taubes:99}.
This permuting action descends to the \hK quotient $\cM_0$ so that Fueter sections with values in $\cM_0$ also make sense on any three--manifold.
In this case, however, one should consider sections of an appropriate bundle rather than just maps into $\cM_0$. 
More details on this can be found in~\cite{Haydys_ahol:08}.
\end{remark}

\section{Compactness of the moduli space of the Seiberg--Witten monopoles with multiple spinors}\label{Sect_CompactnessSW}

In this section $\cM$ is always assumed to be a vector space so that there is no need to restrict to the flat space as the background three-manifold. 
Thus, assume that $M$ is an arbitrary closed oriented Riemannian three-manifold. 
Fix a $\Spin$--structure $\fs$ on $M$ and denote by $\slS$ the associated spinor bundle;
also fix a $\U(1)$--bundle $\sL$ over $M$, a positive integer $n \in \N$ and an $\SU(n)$--bundle $E$ together with a connection $B$.
Consider pairs $(A,\Psi) \in \sA(\sL) \times \Gamma\(\Hom(E,\slS \otimes \sL)\)$ consisting of a connection $A$ on $\sL$ and an $n$-tuple of twisted spinors $\Psi$ satisfying the \emph{Seiberg--Witten equations with $n$ spinors}:
\begin{equation}
  \label{Eq_nSW0}
    \slD_{A\otimes B} \Psi = 0 \qand\quad   F_A = \mu(\Psi).
\end{equation}
The map $\mu$ here is given by $\Psi\mapsto \Psi\Psi^*-\tfrac 12 |\Psi|^2$, which is identified with a $2$-form on $M$ in the same fashion as in the classical Seiberg--Witten theory. 
Moreover, if $E$ is just the product bundle so that $\Psi$ can be thought of as an $n$-tuple of spinors $(\psi_1,\dots,\psi_n)$, then $\mu(\Psi)=\sum_j\psi_j\psi_j^*-\tfrac 12 |\psi_j|^2$. 
That is each summand in the last expression is the quadratic map appearing in the $n=1$ Seiberg--Witten equations. 

It is quite easy to see that there are just two possibilities for the moduli space of solutions of~\eqref{Eq_nSW0} to be non-compact: Either there is a sequence of solutions $(A_i, \Psi_i)$, which converges to a reducible solution or the sequence $\| \Psi_i\|_{L^2}$ diverges to $+\infty$.
With this in mind it is natural to blow-up \eqref{Eq_nSW0}, that is, to consider triples $(A,\Psi,\alpha) \in \sA(\sL) \times \Gamma\(\Hom(E,\slS\otimes \sL)\) \times [0,\pi/2]$ satisfying 
\begin{equation}
  \label{Eq_nSW}
  \begin{split}
    \|\Psi\|_{L^2} &= 1, \\
    \slD_{A\otimes B} \Psi &= 0, \qand \\
    \sin(\alpha)^2 F_A &= \cos(\alpha)^2 \mu(\Psi).
  \end{split}  
\end{equation}
The analytically difficult case is when $\e=\tan\alpha \to 0$, since for $\alpha=0$ Equations~\eqref{Eq_nSW} are no longer elliptic, albeit 
equivalent to an elliptic equation.
\begin{thm}[\cite{HaydysWalpuski15_CompThm_GAFA}*{Thm~1.5}]
  \label{Thm_A}
  Let $(A_i,\Psi_i,\alpha_i)$ be a sequence of solutions of \eqref{Eq_nSW}.
  If $\limsup \alpha_i > 0$, then after passing to a subsequence and up to gauge transformations $(A_i,\Psi_i,\alpha_i)$ converges smoothly to a limit $(A,\Psi,\alpha)$.
  If $\limsup \alpha_i = 0$, then after passing to a subsequence the following holds:
  \begin{itemize}
    \item
      There is a closed nowhere-dense subset $Z \subset M$, a connection $A$ on $\sL|_{M\setminus Z}$ and $\Psi \in \Gamma\(M\setminus Z, \Hom(E,\slS\otimes \sL)\)$ such that $(A, \Psi, 0)$ solves \eqref{Eq_nSW}.
      $|\Psi|$ extends to a H\"older continuous function on all of $M$ and $Z = |\Psi|^{-1}(0)$.
    \item
      On $M \setminus Z$, up to gauge transformations, $A_i$ converges weakly in $W^{1,2}_\loc$ to $A$ and $\Psi_i$ converges weakly in $W^{2,2}_\loc$ to $\Psi$.
      There is a constant $\gamma > 0$ such that $|\Psi_i|$ converges to $|\Psi|$ in $C^{0,\gamma}$ on all of $M$.\qed
  \end{itemize}
\end{thm}

\begin{rem}
Very recently Taubes~\cite{Taubes16_SWDim4_Arx} proved an analogue of Theorem~\ref{Thm_A} in dimension four.
\end{rem}

It follows from the proof of Theorem~\ref{Thm_A} that outside of $Z$ the local $L^2$-energy of $F_A$, i.e., $\int_{B_r(m)}|F_A|^2$, remains bounded; Here $B_r(m)$ is the geodesic ball of radius $r$ centered at $m\in M\setminus Z$. 
Hence, I call $Z$ \emph{a blow-up set} for~\eqref{Eq_nSW} (although the set where the local $L^2$-energy actually blows up may be a proper subset of $Z$).

Let $\mathring M_{1,n}$ be the framed moduli space of centered charge one $\SU(n)$ ASD instantons on $\R^4$.
By the ADHM-construction of instantons on $\R^4$ \cite{DonaldsonKronheimer:90}*{Sect.\,3.3} the \hK  quotient $\H^n\setminus\{ 0 \}\hkred \rU(1)$ is diffeomorphic (in fact, isometric) to $\mathring M_{1,n}$.  
Then $\H^n\hkred \rU(1)$ is the Uhlenbeck compactification $\bar M_{1,n}$ of $\mathring M_{1,n}$ and the image of the origin is a singular point in $\bar M_{1,n}$.
For example, in the case $n=2$ we have $\bar M_{1,2}=\H/\pm 1.$  
 
It follows from (a slight modification of) Theorem~\ref{Thm_FuetMapsHKred} that  the limit $(A,\Psi,0)$ appearing in Theorem~\ref{Thm_A} corresponds to a Fueter section of a bundle $\fM$ with fiber $\mathring M_{1,n}$.
More precisely, this means that  $(A,\Psi,0)$ yields a pair $(\psi, Z)$ such that
\begin{itemize}[itemsep=-2pt,topsep=2pt]
\item $Z$ is a closed nowhere  dense subset of $M$;
\item $\psi$ is a Fueter section of $\fM$ over $M\setminus Z$;
\item $|\psi|$ extends as a continuous function on all of $M$ and $Z=|\psi|^{-1}(0)$.
\end{itemize}
Here $|\cdot|\colon \mathring M_{1,n}\to \R$ denotes the function induced by the Euclidean norm on $\H^n$.
For example, for $n=2$ pairs $(\psi, Z)$ constitute $\Z/2$--harmonic spinors studied in~\cite{Taubes14_ZeroLoci_Arx}.

\medskip

Let me consider the case $n=2$ in some detail.
Since $M$ is three-dimensional and the singularity in $\H/\pm 1$ is of codimension four one may na\"{i}vely hope that $Z$ is empty at least generically.  
The following argument shows that this is not the case.
Indeed, let $(A,\Psi, 0)$ be a solution of~\eqref{Eq_nSW} with $n=2$.
Assume $Z$ is empty. 
The equation $\Psi\Psi^*=\tfrac 12 |\Psi|^2\neq 0$ then implies that $\ker\Psi^*$ is trivial, i.e., $\Psi$ is an epimorphism everywhere.
By dimensional reasons, $\Psi$ must be an isomorphism.  
Hence, $\Lambda^2 E\cong \Lambda^2(\slS\otimes \sL)\cong \sL^2$, which implies that the determinant line bundle $\sL^2$ must be trivial, since $E$ is an $\SU(2)$--bundle. 
Thus, if the determinant line bundle is non-trivial, then $Z$ can not be empty.

More generally, it can be shown that $Z$ carries an extra structure $(\theta, or)$, which consists of the multiplicity function and an orientation in the case of smooth $Z$. 
This extra structure can be used to define the homology class $[Z,\theta, or]\in H_1(M,\Z)$.
More details can be found in~\cite{Haydys16_TopBlowUpSet_Arx}.
\begin{thm}[\cite{Haydys16_TopBlowUpSet_Arx}]
  \label{Thm_TopologyZ}
If $(A,\Psi,0)$ is a solution of~\eqref{Eq_nSW} with $n=2$ over $M\setminus Z$ such that $Z=|\Psi|^{-1}(0)$, then
\begin{equation}
\tag*{$\square$}
[Z,\theta, or]=\PD (c_1(\sL^2)).
\end{equation}
\end{thm}

This theorem immediately implies that there are restrictions on $\Z/2$-harmonic spinors  which can potentially appear as degenerations of the Seiberg--Witten monopoles with $n=2$ (in the sense of Theorem~\ref{Thm_A}) for a fixed choice of the determinant line bundle.
For instance, an honest harmonic spinor with a finite zero locus, viewed as a $\Z/2$-harmonic spinor,  never appears as a degeneration of the Seiberg--Witten monopoles with $n=2$ provided the determinant line bundle is non-trivial.

\section{Compactness of the moduli space of $\rG_2$ instantons}\label{Sect_CompG2Instant}

In this section I describe the compactness problem for the moduli space of $\rG_2$ instantons as well as the r\^{o}le of Fueter sections in this problem.

Pick a compact Lie group $\rK$ and a principal $\rK$-bundle $P\to Y$. 
Recall that a connection $A$ on $P$ is said to be a $\rG_2$ instanton, if $\pi_7(F_A)=0$.  
As already mentioned in the introduction, $\rG_2$ instantons have a lot in common with flat connections on three-manifolds. 
However, while the moduli space of flat connections with a compact gauge group on a closed three-manifold is always compact, the moduli space of $\rG_2$ instantons does not need to have this property. 
This is explained in some details below.

\medskip

Let $A_i$ be a sequence of $\rG_2$ instantons on $P$ over a compact $\rG_2$ manifold $Y$. 
Then~\cites{Uhlenbeck82_ConnWithBounds, Price83_MonotonicityFormula, Nakajima88_CompHighDim} there exists a subset $M\subset Y$ of dimension at most $3$ such that a subsequence of $A_i$ converges over $Y\setminus M$ to a $\rG_2$ instanton (after applying gauge transformations if necessary). 
In general, the limiting instanton may have unremovable singularities over a subset of $M$ and\,/\,or the set $M$ may have singularities.
However, for the purposes of this discussion I assume that $M$ is a smooth submanifold of dimension three and the limiting $\rG_2$ instanton has no unremovable singularities.
Then~\cite{Tian:00} $M$ is \emph{an associative submanifold} of $Y$, i.e., 
\begin{equation}\label{Eq_AssocCond}
\imath_M^*\varphi=\mathrm{vol_M},
\end{equation}
where $\mathrm{vol_M}$ is the induced volume form on $M$. 
Notice that the normal bundle $N_M$ is isomorphic to the (twisted) spinor bundle of $M$~\cite{McLean:98}. 
In particular, there is a Clifford multiplication $T^*M\otimes N_M\to N_M$ and, hence, a Dirac operator $\dirac\colon\Gamma(N_M)\to\Gamma(N_M)$, whose kernel is the space of infinitesimal deformations of $M$ as an associative submanifold.
Hence, generically associative submanifolds should be isolated.

Assume for the sake of simplicity of exposition that $Y=\mathbb R^7$, which is considered as a $\rG_2$ manifold\footnote{Of course, $\mathbb R^7$ is not compact so that strictly speaking the above discussion does not apply. However, it is possible to proceed with any compact $Y$ at the cost of technical complications, which are immaterial for the purposes of this discussion.} with its Euclidean metric $g_{0}$, and $M=\mathbb R^3\times\{ 0 \}$. 
Then a connection $A$ on the product $\rK$-bundle can be thought of as a 1-form on $\mathbb R^7=\mathbb R^3\times \mathbb R^4$ and therefore can be split into the K\"unneth-type components: $A=a+b$. 
Here $a\in \Omega^{1,0}(\mathbb R^7;\fk)$ vanishes on tangent vectors to $\mathbb R^4$, $b\in\Omega^{0,1}(\mathbb R^7;\fk)$ vanishes on tangent vectors to $\mathbb R^3$ and $\fk=\mathrm{Lie}(\rK)$.  

Think of $b$ as a family of 1-forms on $\mathbb R^4$ by identifying $\Omega^{0,1}(\mathbb R^7;\fk)$ with $C^\infty (\mathbb R^3; \Omega^1(\mathbb R^4;\fk))$ and denote by $F_b$ the corresponding family of curvatures.
To be more precise, we are interested only in those connections $A$, which decay sufficiently fast along normal fibres. 
In particular, $b$ takes values in  a subspace $\hat\Omega^1(\mathbb R^4;\fk)$ of 1-forms decaying at $\infty$ (see~\cite{Itoh:89} for details).

Using the identification $\Omega^{1,0}(\mathbb R^7;\fk)\cong\Omega^1(\mathbb R^3; \hat\Omega^0(\mathbb R^4;\fk))$ and observing that $\hat\Omega^0(\mathbb R^4;\fk)$ is the Lie algebra of the based gauge group $\hat\cG\subset\mathrm{Aut}(\underline{\rK}\to\mathbb R^4)$, the component $a$ can be thought of as a connection on the trivial bundle $\underline{\hat\cG}\to\mathbb R^4$ with an infinite dimensional structure group.
Then a computation yields that $A$ is a $\rG_2$ instanton with respect to $g_0$ if and only if  
\begin{equation}\label{Eq_G2InstantOnR7}
I_1\nabla_1^a b + I_2\nabla_2^a b+ I_3\nabla_3^a b=0\quad\text{and}\quad
F_a=F_b^+.
\end{equation}
Here $I_j$'s are complex structures on $\hat\Omega^1(\mathbb R^4;\fk)$ inherited from $\mathbb R^4$ and in the second equation the identification $\Omega^2(\mathbb R^3;\Omega^0(\mathbb R^4;\fk))\cong C^\infty(\mathbb R^7;\fk\otimes\mathbb R^3)\cong C^\infty (\mathbb R^3;\Omega^2_+(\mathbb R^4;\fk))$ is used. 

For any $\e>0$ denote by $\rho_\e\colon \mathbb R^7\to \mathbb R^7$ the dilation along the normal fibres: $(x,y)\mapsto (x,\e y)$. 
Clearly, if $A$ is a $\rG_2$ instanton with respect to $g_0$, then $A_\e:=\rho_\e^*A$ is a $\rG_2$ instanton with respect to 
$g^\e:=\rho_\e^*g_0=dx^2+\e^2dy^2$.
Arguing along similar lines as above, we obtain that with $A_\e =a_\e+b_\e$ this is equivalent to 
\begin{equation}\label{Eq_G2InstantOnR7e}
 I_1\nabla_1^{a_\e} b_\e + I_2\nabla_2^{a_\e} b_\e+ I_3\nabla_3^{a_\e} b_\e=0\quad\text{and}\quad
 \e^2F_{a_\e}=F_{b_\e}^+.
\end{equation}  
Clearly, if for some sequence of $\rG_2$ instantons $A_i$ developing a bubble along $\mathbb R^3$ the rescaled sequence $\rho_{\e_i}^*A_i$ converges for a suitably chosen $\e_i\to 0$, then the limit $A_0=a_0+ b_0$ is a solution of
\begin{equation}\label{Eq_G2InstantOnR70}
 I_1\nabla_1^{a_0} b_0 + I_2\nabla_2^{a_0} b_0+ I_3\nabla_3^{a_0} b_0=0\quad\text{and}\quad
 F_{b_0}^+=0.
\end{equation}  

To understand the meaning of~\eqref{Eq_G2InstantOnR70}, observe that $\hat\Om^1(\mathbb R^4;\fk)$ can be viewed as an infinite dimensional hyperK\"ahler manifold. 
Moreover, the action of the based gauge group preserves the hyperK\"ahler structure and the corresponding moment map is given by $b\mapsto F_b^+$.
Hence, \eqref{Eq_G2InstantOnR70} is a special case\footnote{It is also interesting to observe that~\eqref{Eq_G2InstantOnR7} and~\eqref{Eq_G2InstantOnR7e} are special cases of~\eqref{Eq_gSW} and~\eqref{Eq_gSW_e} respectively.} of~\eqref{Eq_gSW_0} with $(\cM,\rG)=(\hat\Om^1(\R^4;\fk),\hat\cG)$.
It is easy to see that Theorem~\ref{Thm_FuetMapsHKred} also holds in this infinite dimensional case and, hence, if the limit of the rescaled sequence of $\rG_2$ instantons developing a bubble along $\mathbb R^3$ exists, this corresponds to a Fueter map with values in the framed moduli space of anti-self-dual instantons over $\mathbb R^4$.  

The whole discussion can be generalized for an arbitrary smooth associative (compact) submanifold $M$ of a $\rG_2$ manifold $Y$ (details are spelled out in~\cite{Haydys:12_GaugeTheory_jlms} for the case of $\Spin(7)$-manifolds).  
Namely, replace each fiber $N_mM\cong\mathbb R^4$ of the normal fiber bundle by the framed moduli space of anti-self-dual instantons over $N_mM$ to obtain a fiber bundle $\fM\to M$.
There is an intrinsic Fueter equation for sections of $\fM$ and if appropriately rescaled sequence of $\rG_2$ instantons, which develops a bubble along $M$, converges, then the limit is a Fueter section.

\medskip

At present there is no rigorous proof that the bubbles of $\rG_2$ instantons are indeed modelled on Fueter sections. 
However, Walpuski~\cite{Walpuski12_G2AndFueter_arx} showed that a generic (in a suitable sense) Fueter section on an associative submanifold of a compact $\rG_2$ manifold can be deformed into a family of $\rG_2$ instantons, whose bubble is modeled on the initial Fueter section.

\section{A relation between $\rG_2$ instantons and generalized Seiberg--\-Wit\-ten equations}\label{Sect_RelationG2SW}

The material of this section is mostly speculative. 
The main purpose is to describe the qualitative picture of the conjectural relation between  $\rG_2$ instantons, Fueter sections, and generalized Seiberg--Witten monopoles. 
A quantitative version of this relation is an ongoing project with Walpuski; Details will be described elsewhere.

Let $\cM_{asd}$ denote the framed moduli space of all anti-self-dual instantons on $S^4=\R^4\cup\{\infty\}$.
The action of $\R_{>0}$ on $\R^4$ by dilations induces an $\R_{>0}$-action on $\cM_{asd}$.
Since this action preserves all complex structures, the Fueter sections with values in $\cM_{asd}$ always come in one-parameter families. 
Hence, given a one-parameter family of metrics on a fixed three-manifold $M$ we expect that Fueter sections with values in $\cM_{asd}$ appear generically for a finite number of metrics from this family. 

\begin{figure}[ht]
   \begin{center}
    \includegraphics[width=0.75\textwidth]{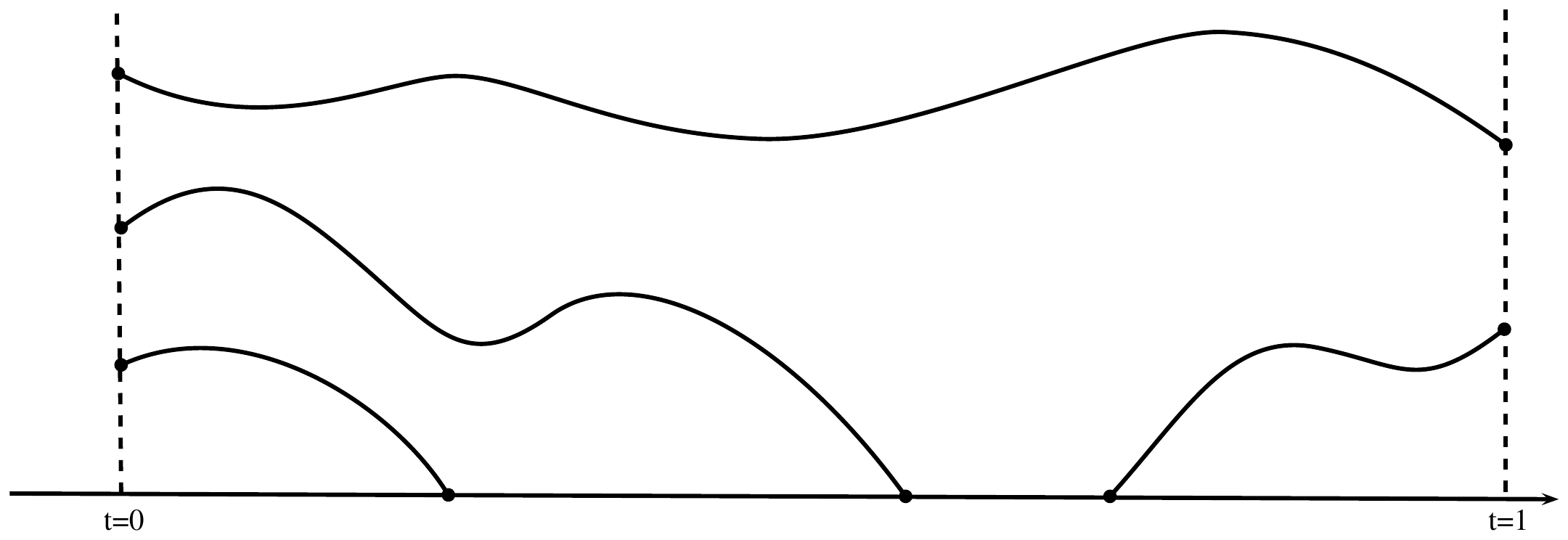}
   \end{center}
\caption{Expected behaviour of $\rG_2$ instantons in one-parameter families.}
\label{Fig_G2instInFam}
\end{figure}

The above means that for a generic $\rG_2$ metric (in a suitable sense) the moduli space of $\rG_2$ instantons is expected to consist of finitely many points. 
However in a one-parameter family $g_t$ of such metrics we expect to have finitely many metrics, for which one of the instantons may develop a bubble, become a Fueter section in the limit and die afterwards.
Symmetrically, some $\rG_2$ instantons can be born out of a Fueter section.
This is schematically depicted on Fig.\ref{Fig_G2instInFam}.

The upshot of this discussion is that the algebraic number of $\rG_2$ instantons is not expected to be invariant with respect to deformations of the underlying $\rG_2$ structure.  
This was described in~\cite{DonaldsonSegal:09}*{Sect.\,6.2} for the first time and a question about possible counterterm compensating vanishings and appearances of $\rG_2$ instantons was raised. 
Below I sketch an approach to the construction of such counterterm, which is based on generalized \SWn equations rather than a non-linear spectral flow as originally suggested by Donaldson and Segal.  

\medskip

\begin{figure}[ht]
   \begin{center}
    \includegraphics[width=0.75\textwidth]{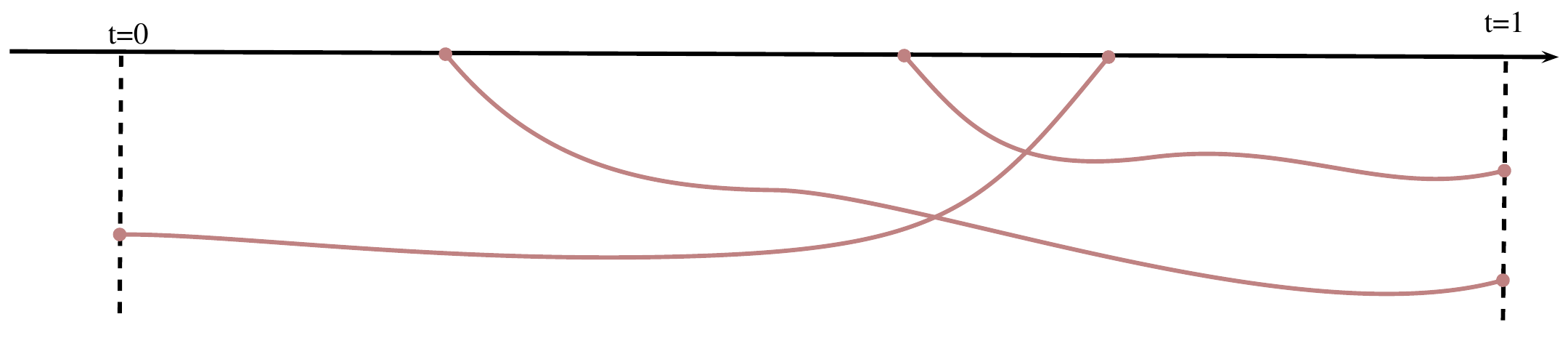}
   \end{center}
\caption{Expected behaviour of Seiberg--Witten monopoles in one-parameter families.}
\label{Fig_SWInFam}
\end{figure}

A relation between $\rG_2$ instantons and generalized \SWn equations is based on the observation that (conjecturally)  the behavior of generalized Seiberg--Witten monopoles on a fixed three-manifold in a one-parameter family resembles Fig.\ref{Fig_G2instInFam}. 
More precisely, in the setting of Section~\ref{Sect_CompactnessSW} for a generic metric on the background three-manifold $M$ the moduli space of solutions of~\eqref{Eq_nSW0} should consist of finitely many isolated points. 
However, in a one-parameter family there should be a finite number of metrics, for which a generalized Seiberg--Witten monopole converges to a Fueter section and disappears afterwards or, symmetrically, a generalized \SWn monopole can be born out of a Fueter section. This is schematically depicted on Fig.\ref{Fig_SWInFam}.

Putting these two pictures together, the following interplay between $\rG_2$ instantons and generalized Seiberg--Witten monopoles is expected. 
Let $g_t,  \ t\in [0,1]$, be a one-parameter family of $\rG_2$ metrics on a compact manifold $Y$. 
Assume $A_t, t\in [0, t_0)$, is a one-parameter family of $\rG_2$ instantons, which develops a bubble at $t=t_0$ along an associative submanifold $\bigl (M, \left. g_{t_0}\right |_{M} \bigr )\subset (Y, g_{t_0})$ and converges to a Fueter section on $M$ with values in $\cM_{asd}$.   
This Fueter section gives birth to a one-parameter family of generalized \SWn monopoles on $M$ for $t>t_0$ (see Fig.\ref{Fig_G2andSWinstInFam} below).
Of course, we should also have a symmetric process, in which a family of generalized \SWn monopoles degenerates to a Fueter section and gives  birth to $\rG_2$ instantons. 
\begin{figure}[ht]
   \begin{center}
    \includegraphics[width=0.75\textwidth]{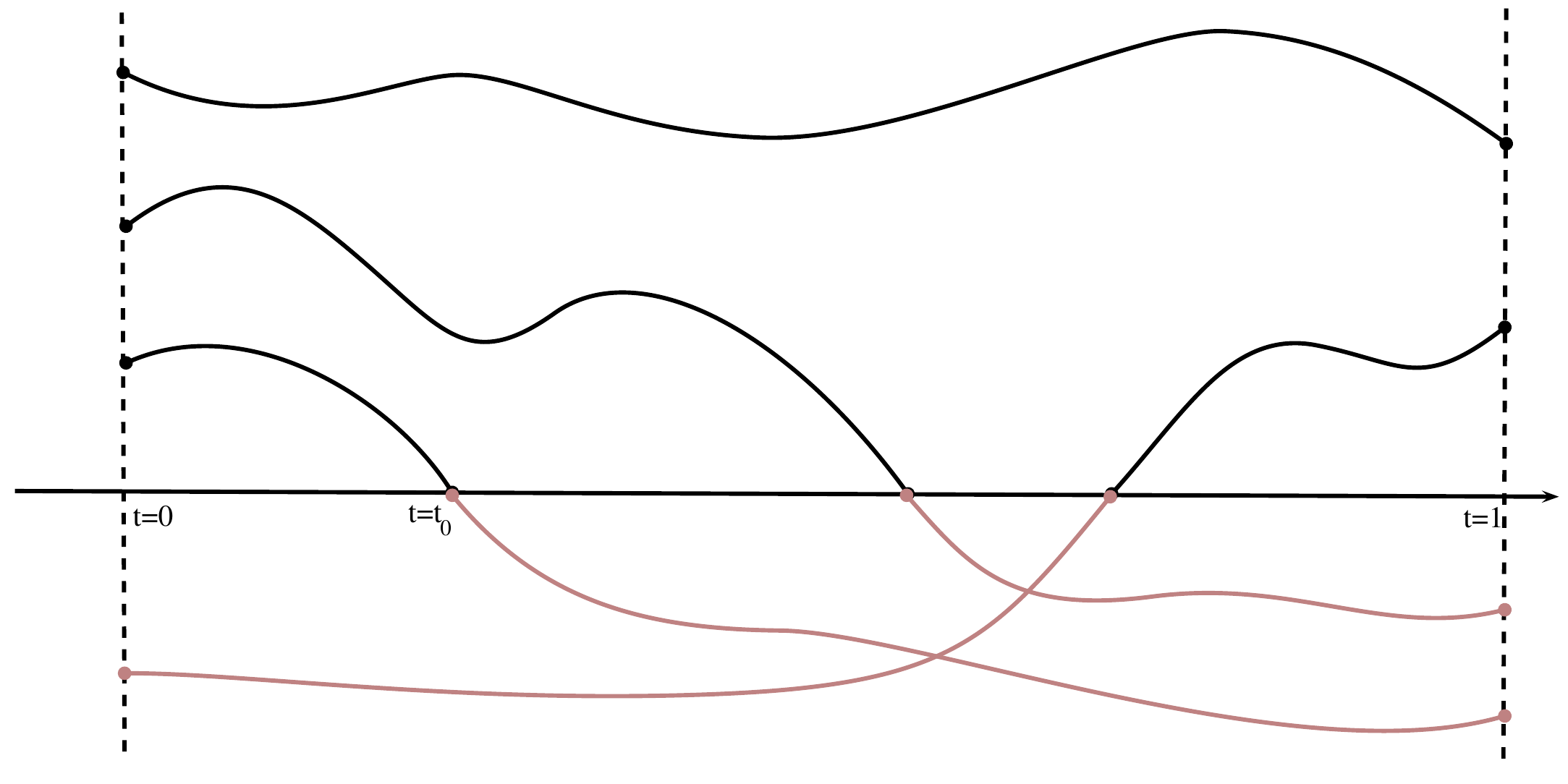}
   \end{center}
\caption{Expected behavior of $\rG_2$ instantons and Seiberg--Witten monopoles on associative submanifolds in one-parameter families.}
\label{Fig_G2andSWinstInFam}
\end{figure}

This suggests that in order to construct an invariant out of $\rG_2$ instantons, which would remain constant along one-parameter families of $\rG_2$ metrics, one has to take generalized \SWn monopoles on associative submanifolds into account. 
More precisely, for a bundle $E\to Y$ such an invariant (if exists) should be schematically of the form
\begin{equation}
  \label{Eq_G2GassonInvariant}
  n(E) := \sum_{([A],M)} w([A],M,g),
\end{equation}   
where $[A]$ is the gauge equivalence class of $\rG_2$ instantons, $M\subset Y$ is an associative submanifold and $w$ is a weight.
In particular, $w([A],\varnothing ,g)=\pm 1$ and for $M\neq\varnothing$ the corresponding weight should be constructed by counting certain generalized \SWn monopoles. 
The simplest example of such monopoles is~\eqref{Eq_nSW0}.
%This is a joint project with Walpuski; Details will appear elsewhere.

\medskip

It is clear that the picture outlined above may be very rough approximation to reality because, for example, $M$ does not need to be smooth in general or in a one-parameter family of $\rG_2$ metrics corresponding associative submanifolds can undergo a geometric transition and apriori it is not clear whether~\eqref{Eq_G2GassonInvariant} remains invariant under such transitions.  
Recently Joyce~\cite{Joyce16_ConjCountingAssoc_Arx} pointed out that it is not even clear whether a counter-term to the count of $\rG_2$ instantons can be defined in terms of counting associative submanifolds with some weights. 
While it is not clear indeed, the argument in \cite{Joyce16_ConjCountingAssoc_Arx} uses implicitly a hypothesis that (in the simplest case) $\Z/2$-harmonic spinors behave like honest harmonic spinors.
The authors feeling is that this may not be the case.
In any case the question to which extend the known facts for honest harmonic spinors apply also to $\Z/2$-harmonic spinors requires further studies.

\section{Some open questions} 

The purpose of this section is to gather some open questions related to the content of the  preceding sections.

One broad question related to Theorem~\ref{Thm_A} is the following: Which pairs $(\psi, Z)$ can be realized as degenerations of the Seiberg--Witten monopoles with two spinors (more generally, with $n$ spinors)? 
As already mentioned at the end of Section~\ref{Sect_CompactnessSW}, not every harmonic $\Z/2$-spinor can appear as a degeneration of the Seiberg--Witten monopoles with two spinors at least if the determinant line bundle is fixed.
The above question includes in particular the following.
\begin{question}
What can be said about the regularity of $Z$?
\end{question}

If $\psi$ is an honest harmonic spinor, then the nodal set $\psi^{-1}(0)$ is known to be countably rectifiable \cite{Baer97_NodalSetsDirac} in the sense that it is contained in the countable union of the images of Lipschitz maps $\R\to M$.
In particular, the Hausdorff dimension of the nodal set is at most one.
For $\Z/2$-harmonic spinors a known variant of this \cite{Taubes14_ZeroLoci_Arx} is that $Z$ contains an open everywhere dense subset, which is a Lipschitz submanifold of $M$. 
Taubes also proves that the Hausdorff dimension of $Z$ is at most one. 
However the question whether $Z$ is rectifiable remains open. It is also not known whether the one dimensional Hausdorff measure of $Z$ is finite.  

\medskip

If the determinant line bundle is trivial then the argument presented at the end of Section~\ref{Sect_CompactnessSW} does not apply and we can ask the following question.
\begin{question}
Let $(A,\Psi,0)$ be a solution of~\eqref{Eq_nSW} corresponding to a trivial  determinant line bundle. 
Is the blow-up set $Z$ empty generically?
\end{question}

Generic in this question means for an open dense subset in the space $\cR\times\cA(E)$, where $\cR$ is the space of all Riemannian metrics on $M$ and $\cA(E)$ is the space of all $\SU(n)$-connections on $E$.

This question seems to be particularly relevant for the relation with the $\rG_2$ instantons.
Indeed, assume a sequence of $\rG_2$ instantons developing a bubble along an associative submanifold converges to a Fueter section $\psi$ with values in $\bar M_{1,2}=\H/\pm 1$, i.e., $\psi$ is a $\Z/2$-harmonic spinor.
Then $\psi$ may hit the singularity in $\H/\pm 1$ at some points but there are no reasons to expect that the preimage of the singular point determines a  non-trivial homology class.    
Hence, it seems reasonable to expect that only the solutions $(A,\Psi, 0)$ with the trivial class $[Z,\theta, or]$ are relevant for the relation with $\rG_2$ instantons. 
If $Z$ could be assumed to be empty generically (both in the case of the Seiberg--Witten monopoles and $\rG_2$ instantons), this would simplify a lot technical work needed to be done to prove the relation.

\medskip

If $n>2$, then a variant of Theorem~\ref{Thm_TopologyZ} can be obtained as follows. 
Let $(A,\Psi,0)$ be a solution of~\eqref{Eq_nSW} over $M\setminus Z$.
It is convenient to choose a topological trivialization of $E$ thus identifying $E$ with the product bundle $\underline\C^n$ (a trivialization exists since the base manifold is three-dimensional).
The equation $\mu(\Psi)=0$ implies that $\Psi$ is an epimorphism over $M\setminus Z$. 
Thus we obtain a map
\[
\Phi_0\colon M\setminus Z\to \Gr_{n-2}(\C^n),\qquad \Phi_0(m)=\ker\Psi_m.
\]
If $\Phi_0$ admits an extension to a map defined on all of $M$, then one can still define $(\theta, or)$ on $Z$ so that the class $[Z,\theta, or]\in H_1(M,\Z)$ is well-defined. 
\begin{question}
Does the map $\Phi_0$ extend to a continuous map on all of $M$ for any solution $(A,\Psi,0)$ of~\eqref{Eq_nSW} over $M\setminus Z$?
\end{question}

\bibliography{references}

\def\cprime{$'$}
\begin{bibdiv}
\begin{biblist}

\bib{Baer97_NodalSetsDirac}{article}{
      author={B{\"a}r, C.},
       title={On nodal sets for {D}irac and {L}aplace operators},
        date={1997},
        ISSN={0010-3616},
     journal={Comm. Math. Phys.},
      volume={188},
      number={3},
       pages={709\ndash 721},
         url={http://dx.doi.org/10.1007/s002200050184},
      review={\MR{1473317}},
}

\bib{Berger55_HolGps}{article}{
      author={Berger, M.},
       title={Sur les groupes d'holonomie homog\`ene des vari\'et\'es \`a
  connexion affine et des vari\'et\'es riemanniennes},
        date={1955},
        ISSN={0037-9484},
     journal={Bull. Soc. Math. France},
      volume={83},
       pages={279\ndash 330},
      review={\MR{0079806 (18,149a)}},
}

\bib{Bryant87_MetrWithExceptHol}{article}{
      author={Bryant, R.},
       title={Metrics with exceptional holonomy},
        date={1987},
        ISSN={0003-486X},
     journal={Ann. of Math. (2)},
      volume={126},
      number={3},
       pages={525\ndash 576},
         url={http://dx.doi.org/10.2307/1971360},
      review={\MR{916718 (89b:53084)}},
}

\bib{BryantSalamon89_CompleteExcept}{article}{
      author={Bryant, R.},
      author={Salamon, S.},
       title={On the construction of some complete metrics with exceptional
  holonomy},
        date={1989},
        ISSN={0012-7094},
     journal={Duke Math. J.},
      volume={58},
      number={3},
       pages={829\ndash 850},
         url={http://dx.doi.org/10.1215/S0012-7094-89-05839-0},
      review={\MR{1016448 (90i:53055)}},
}

\bib{SalamonEtAl:02}{article}{
      author={Cieliebak, Kai},
      author={Gaio, A.~Rita},
      author={{Mundet i Riera}, Ignasi},
      author={Salamon, Dietmar~A.},
       title={{The symplectic vortex equations and invariants of {H}amiltonian
  group actions}},
        date={2002},
     journal={J. Symplectic Geom.},
      volume={1},
      number={3},
       pages={543\ndash 645},
}

\bib{CGN_AnalytInvOfG2mflds_arx}{article}{
      author={Crowley, D.},
      author={Goette, S.},
      author={Nordstr{\"o}m, J.},
       title={An analytic invariant of $g_2$ manifolds},
        date={2015},
     journal={arXiv:1505.02734},
}

\bib{CHNP15_G2Mflds}{article}{
      author={Corti, A.},
      author={Haskins, M.},
      author={Nordstr{\"o}m, J.},
      author={Pacini, T.},
       title={{$\rm{G}_2$}-manifolds and associative submanifolds via
  semi-{F}ano 3-folds},
        date={2015},
        ISSN={0012-7094},
     journal={Duke Math. J.},
      volume={164},
      number={10},
       pages={1971\ndash 2092},
         url={http://dx.doi.org/10.1215/00127094-3120743},
      review={\MR{3369307}},
}

\bib{DonaldsonKronheimer:90}{book}{
      author={Donaldson, S.~K.},
      author={Kronheimer, P.~B.},
       title={{The geometry of four-manifolds}},
      series={{Oxford Mathematical Monographs}},
   publisher={The Clarendon Press Oxford University Press},
     address={New York},
        date={1990},
        ISBN={0-19-853553-8},
        note={Oxford Science Publications},
}

\bib{DonaldsonSegal:09}{incollection}{
      author={Donaldson, S.},
      author={Segal, E.},
       title={Gauge theory in higher dimensions, {II}},
        date={2011},
   booktitle={Surveys in differential geometry. {V}olume {XVI}. {G}eometry of
  special holonomy and related topics},
      series={Surv. Differ. Geom.},
      volume={16},
   publisher={Int. Press, Somerville, MA},
       pages={1\ndash 41},
      review={\MR{2893675}},
}

\bib{DonaldsonThomas:98}{incollection}{
      author={Donaldson, S.~K.},
      author={Thomas, R.~P.},
       title={{Gauge theory in higher dimensions}},
        date={1998},
   booktitle={{The geometric universe (Oxford, 1996)}},
   publisher={Oxford Univ. Press},
     address={Oxford},
       pages={31\ndash 47},
}

\bib{Haydys_ahol:08}{article}{
      author={Haydys, A.},
       title={{Nonlinear Dirac Operator and Quaternionic Analysis}},
        date={2008},
        ISSN={0010-3616},
     journal={Communications in Mathematical Physics},
      volume={281},
      number={1},
       pages={251\ndash 261},
         url={http://dx.doi.org/10.1007/s00220-008-0445-1},
}

\bib{Haydys:12_GaugeTheory_jlms}{article}{
      author={Haydys, A.},
       title={Gauge theory, calibrated geometry and harmonic spinors},
        date={2012},
        ISSN={0024-6107},
     journal={J. Lond. Math. Soc. (2)},
      volume={86},
      number={2},
       pages={482\ndash 498},
         url={http://dx.doi.org/10.1112/jlms/jds008},
      review={\MR{2980921}},
}

\bib{Haydys15_DiracOperatorsInGauge}{incollection}{
      author={Haydys, A.},
       title={Dirac operators in gauge theory},
        date={2015},
   booktitle={New ideas in low dimensional topology},
      series={Ser. Knots Everything},
      volume={56},
   publisher={World Sci. Publ., Hackensack, NJ},
       pages={161\ndash 188},
         url={http://dx.doi.org/10.1142/9789814630627_0005},
      review={\MR{3381325}},
}

\bib{Haydys16_TopBlowUpSet_Arx}{article}{
      author={Haydys, A.},
       title={Topology of the blow-up set for the {S}eiberg--{W}itten equation
  with multiple spinors},
        date={2016},
      eprint={\arxiv{1607.01763}},
}

\bib{HaydysWalpuski15_CompThm_GAFA}{article}{
      author={Haydys, A.},
      author={Walpuski, Th.},
       title={A compactness theorem for the {S}eiberg--{W}itten equation with
  multiple spinors in dimension three},
        date={2015},
        ISSN={1016-443X},
     journal={Geom. Funct. Anal.},
      volume={25},
      number={6},
       pages={1799\ndash 1821},
         url={http://dx.doi.org/10.1007/s00039-015-0346-3},
        note={Erratum available via
  \href{http://dx.doi.org/10.13140/RG.2.1.2559.1285 }{doi:
  10.13140/RG.2.1.2559.1285}},
      review={\MR{3432158}},
}

\bib{Itoh:89}{incollection}{
      author={Itoh, M.},
       title={Based anti-instantons and gravitational instantons},
        date={1989},
   booktitle={Geometry of manifolds ({M}atsumoto, 1988)},
      series={Perspect. Math.},
      volume={8},
   publisher={Academic Press},
     address={Boston, MA},
       pages={453\ndash 475},
}

\bib{Joyce07_RiemannianHolGpsCalibr}{book}{
      author={Joyce, D.},
       title={Riemannian holonomy groups and calibrated geometry},
      series={Oxford Graduate Texts in Mathematics},
   publisher={Oxford University Press, Oxford},
        date={2007},
      volume={12},
        ISBN={978-0-19-921559-1},
      review={\MR{2292510}},
}

\bib{Joyce16_ConjCountingAssoc_Arx}{article}{
      author={Joyce, D.},
       title={Conjectures on counting associative 3-folds in
  {$\rG_2$}-manifolds},
        date={2016},
      eprint={\arxiv{1610.09836}},
}

\bib{Joyce96_CompactG2I_II}{article}{
      author={Joyce, D.},
       title={Compact {R}iemannian {$7$}-manifolds with holonomy {$G_2$}. {I},
  {II}},
        date={1996},
        ISSN={0022-040X},
     journal={J. Differential Geom.},
      volume={43},
      number={2},
       pages={291\ndash 328, 329\ndash 375},
         url={http://projecteuclid.org/euclid.jdg/1214458109},
      review={\MR{1424428 (97m:53084)}},
}

\bib{Joyce96_CompactSpin7}{article}{
      author={Joyce, D.~D.},
       title={Compact {$8$}-manifolds with holonomy {${\rm Spin}(7)$}},
        date={1996},
        ISSN={0020-9910},
     journal={Invent. Math.},
      volume={123},
      number={3},
       pages={507\ndash 552},
         url={http://dx.doi.org/10.1007/s002220050039},
      review={\MR{1383960 (97d:53052)}},
}

\bib{Kovalev03_TCS}{article}{
      author={Kovalev, A.},
       title={Twisted connected sums and special {R}iemannian holonomy},
        date={2003},
        ISSN={0075-4102},
     journal={J. Reine Angew. Math.},
      volume={565},
       pages={125\ndash 160},
         url={http://dx.doi.org/10.1515/crll.2003.097},
      review={\MR{2024648 (2004m:53088)}},
}

\bib{McLean:98}{article}{
      author={McLean, R.},
       title={{Deformations of calibrated submanifolds}},
        date={1998},
        ISSN={1019-8385},
     journal={Comm. Anal. Geom.},
      volume={6},
      number={4},
       pages={705\ndash 747},
}

\bib{Nakajima88_CompHighDim}{article}{
      author={Nakajima, H.},
       title={Compactness of the moduli space of {Y}ang-{M}ills connections in
  higher dimensions},
        date={1988},
        ISSN={0025-5645},
     journal={J. Math. Soc. Japan},
      volume={40},
      number={3},
       pages={383\ndash 392},
         url={http://dx.doi.org/10.2969/jmsj/04030383},
      review={\MR{945342 (89g:58050)}},
}

\bib{Price83_MonotonicityFormula}{article}{
      author={Price, P.},
       title={A monotonicity formula for {Y}ang-{M}ills fields},
        date={1983},
        ISSN={0025-2611},
     journal={Manuscripta Math.},
      volume={43},
      number={2-3},
       pages={131\ndash 166},
         url={http://dx.doi.org/10.1007/BF01165828},
      review={\MR{707042 (84m:58033)}},
}

\bib{Salamon03_TourExceptGeom}{article}{
      author={Salamon, S.},
       title={A tour of exceptional geometry},
        date={2003},
        ISSN={1424-9286},
     journal={Milan J. Math.},
      volume={71},
       pages={59\ndash 94},
         url={http://dx.doi.org/10.1007/s00032-003-0015-0},
      review={\MR{2120916}},
}

\bib{Salamon:89}{book}{
      author={Salamon, S.},
       title={{Riemannian geometry and holonomy groups}},
      series={{Pitman Research Notes in Mathematics Series}},
   publisher={Longman Scientific \& Technical},
     address={Harlow},
        date={1989},
      volume={201},
        ISBN={0-582-01767-X},
      review={\MR{MR1004008 (90g:53058)}},
}

\bib{Taubes13_CxASD_Arx}{article}{
      author={Taubes, C.},
       title={Compactness theorems for {$\mathrm{SL}(2;{\Bbb C})$}
  generalizations of the 4-dimensional anti-self dual equations},
        date={2013},
      eprint={\arxiv{1307.6447}},
}

\bib{Taubes13_PSL2Ccmpt}{article}{
      author={Taubes, C.},
       title={{${\rm PSL}(2;\Bbb C)$} connections on 3-manifolds with {${\rm
  L}^2$} bounds on curvature},
        date={2013},
        ISSN={2168-0930},
     journal={Camb. J. Math.},
      volume={1},
      number={2},
       pages={239\ndash 397},
         url={http://dx.doi.org/10.4310/CJM.2013.v1.n2.a2},
      review={\MR{3272050}},
}

\bib{Taubes14_ZeroLoci_Arx}{article}{
      author={Taubes, C.},
       title={The zero loci of {Z}/2 harmonic spinors in dimension 2, 3 and 4},
        date={2014},
      eprint={\arxiv{1407.6206}},
}

\bib{Taubes15_CorrigendumToPSL2C}{article}{
      author={Taubes, C.},
       title={Corrigendum to ``{$PSL(2;\Bbb C)$} connections on 3-manifolds
  with {$L^2$} bounds on curvature'' [{C}ambridge {J}ournal of {M}athematics
  1(2013) 239--397] [ {MR}3272050]},
        date={2015},
        ISSN={2168-0930},
     journal={Camb. J. Math.},
      volume={3},
      number={4},
       pages={619\ndash 631},
      review={\MR{3435274}},
}

\bib{Taubes16_SWDim4_Arx}{article}{
      author={Taubes, C.},
       title={On the behavior of sequences of solutions to {$\rm U(1)$}
  {S}eiberg--{W}itten systems in dimension 4},
        date={2016},
      eprint={\arxiv{1610.07163}},
}

\bib{Taubes17_VafaWitten_Arx}{article}{
      author={Taubes, C.},
       title={The behavior of sequences of solutions to the vafa-witten
  equations},
        date={2017},
      eprint={\arxiv{1702.04610}},
}

\bib{Taubes:99}{incollection}{
      author={Taubes, C.},
       title={{Nonlinear generalizations of a {$3$}-manifold's {D}irac
  operator}},
        date={1999},
   booktitle={{Trends in mathematical physics (Knoxville, TN, 1998)}},
      series={{AMS/IP Stud. Adv. Math.}},
      volume={13},
   publisher={Amer. Math. Soc.},
     address={Providence, RI},
       pages={475\ndash 486},
      review={\MR{MR1708781 (2001a:58048)}},
}

\bib{Tian:00}{article}{
      author={Tian, G.},
       title={{Gauge theory and calibrated geometry. {I}}},
        date={2000},
        ISSN={0003-486X},
     journal={Ann. of Math. (2)},
      volume={151},
      number={1},
       pages={193\ndash 268},
      review={\MR{MR1745014 (2000m:53074)}},
}

\bib{Uhlenbeck82_ConnWithBounds}{article}{
      author={Uhlenbeck, K.},
       title={Connections with {$L^{p}$} bounds on curvature},
        date={1982},
        ISSN={0010-3616},
     journal={Comm. Math. Phys.},
      volume={83},
      number={1},
       pages={31\ndash 42},
         url={http://projecteuclid.org/euclid.cmp/1103920743},
      review={\MR{648356 (83e:53035)}},
}

\bib{Walpuski12_G2AndFueter_arx}{article}{
      author={Walpuski, Th.},
       title={{$G_2$}-instantons, associative submanifolds and {F}ueter
  sections},
        date={2012},
     journal={arXiv:1205.5350},
}

\bib{Wang89_ParallelSpinors}{article}{
      author={Wang, M.},
       title={Parallel spinors and parallel forms},
        date={1989},
        ISSN={0232-704X},
     journal={Ann. Global Anal. Geom.},
      volume={7},
      number={1},
       pages={59\ndash 68},
         url={http://dx.doi.org/10.1007/BF00137402},
      review={\MR{1029845 (91g:53053)}},
}

\end{biblist}
\end{bibdiv}
\end{document}